\documentclass[12pt]{article}
\usepackage{color}
\definecolor{darkblue}{rgb}{0.00,0.25,0.50}
\usepackage[colorlinks, filecolor=blue, citecolor=darkblue]{hyperref}

\usepackage{amsfonts,amssymb,amsmath,amsthm}
\usepackage{mathrsfs}
\usepackage{url}
\usepackage{enumerate}
\usepackage{cite}
\usepackage[ukrainian, russian, english]{babel}
\usepackage[cp1251]{inputenc}
\begin{document} \selectlanguage{ukrainian}
\thispagestyle{empty}

\title{}

UDC 517.51 \vskip 5mm

\begin{center}
\textbf{\Large Найкращі наближення класів періодичних  функцій багатьох змінних  з обмеженою домінуючою мішаною похідною}
\end{center}

\begin{center}
\textbf{\Large Best approximations for classes of periodic functions \\ of many variables with bounded dominating mixed derivative }
\end{center}

\begin{center}
\large{K.\,V.~Pozharska,  A.\,S.~Romanyuk, S.\,Ya.~Yanchenko~\footnote{E-mail:  pozharska.k@gmail.com, romanyuk@imath.kiev.ua, yan.sergiy@gmail.com}} \\
\emph{\small
The Institute of Mathematics of the NAS of Ukraine, Kyiv} 
\end{center}
\begin{center}
\large{К.\,В. Пожарська, А.\,С. Романюк, C.\,Я.~Янченко}\\
\emph{\small Інститут математики НАН України, Київ}
\end{center}
\vskip0.5cm

\begin{abstract}

Встановлено точні за порядком оцінки наближення класів Соболєва $W^{\boldsymbol{r}}_{p,\boldsymbol{\alpha}}(\mathbb{T}^d)$ періодичних функцій багатьох  змінних з обмеженою домінуючою мішаною похідною. Наближення  здійснюється за допомогою  тригонометричних поліномів  зі спектром у східчастих гіперболічних хрестах, а похибка оцінюється в метриці простору $B_{q,1}(\mathbb{T}^d)$, $1 \leqslant p, q < \infty$.

\vskip 3 mm

We  established exact in order estimates an approximation of the Sobolev classes $W^{\boldsymbol{r}}_{p,\boldsymbol{\alpha}}(\mathbb{T}^d)$ of periodic functions of many variables with a bounded dominating mixed derivative. The approximation is made using trigonometric polynomials with the spectrum in step-hyperbolic crosses, and the error is estimated in the metric of the space $B_{q,1}(\mathbb{T}^d)$, $1 \leqslant p, q < \infty$.

\end{abstract}

\emph{\textbf{Ключові слова та фрази:} класи  Соболєва,  найкраще  наближення, домінуюча  мішана  похідна, східчастий гіперболічний  хрест, суми Фур'є.}

\vskip 1 mm
\emph{\textbf{Keywords:} Sobolev classes, best approximation, dominating mixed derivative, step-hyperbolic cross, Fourier sums.}

\vskip 5 mm

\textbf{1. Вступ.} У роботі досліджуються питання апроксимації класів Соболєва $W^{\boldsymbol{r}}_{p,\boldsymbol{\alpha}}(\mathbb{T}^d)$ періодичних функцій   багатьох змінних з обмеженою домінуючою мішаною похідною у просторі $B_{q,1}(\mathbb{T}^d)$, $1 \leqslant p, q < \infty$, норма в якому є більш ``сильною'' ніж $L_q(\mathbb{T}^d)$-норма. В якості агрегатів  наближення функцій з цих класів використовуються тригонометричні  поліноми з гармоніками зі східчастих гіперболічних хрестів. У випадку $1 < q < \infty$ порядки найкращих наближень зазначеними поліномами реалізуються за наближення функцій їхніми східчасто-гіперболічними сумами Фур'є.

Мотивацією до дослідження апроксимаційних характеристик класів $W^{\boldsymbol{r}}_{p,\boldsymbol{\alpha}}(\mathbb{T}^d)$ у просторі $B_{q,1}(\mathbb{T}^d)$, $1 \leqslant q < \infty$, є та обставина, що в результаті досліджень проведених, зокрема, у роботах \cite{Temlyakov-1990Proc-189, Kashin-Temlyakov-1994MN, Belinskii-1998JAT, RomanyukAV-19UMJ2, RomanyukAV-19UMJ8, RomanyukAV-21UMB, Hembarska-Zaderei-2022, Romanyuk-Yanchenko-2021UMJ, Romanyuk-Yanchenko-2022-6UMJ}, де вивчалися питання апроксимації   класів періодичних функцій багатьох змінних мішаної гладкості (Соболєва $W^{\boldsymbol{r}}_{p,\boldsymbol{\alpha}}(\mathbb{T}^d)$, Нікольського--Бєсова $B^{\boldsymbol{r}}_{p,\theta}(\mathbb{T}^d)$) та деяких їхніх  аналогів у просторі $B_{q,1}(\mathbb{T}^d)$, $q \in \{1,\infty \}$, було   зроблено суттєвий просув в оцінках апроксимаційних  характеристик згаданих класів  функцій  у порівнянні з відомими на той  час відповідними оцінками в $L_q$-просторі. Крім  цього також було виявлено, що у переважній більшості ситуацій в багатовимірному випадку, на  відміну  від  одновимірного, відповідні  апроксимаційні характеристики  в    просторах $B_{q,1}(\mathbb{T}^d)$ і $L_q(\mathbb{T}^d)$, $q\in \{1, \infty\}$, мають різні порядки. Таким  чином, з огляду на  сказане, природним є інтерес  до   дослідження апроксимаційних  характеристик класів  $W^{\boldsymbol{r}}_{p,\boldsymbol{\alpha}}(\mathbb{T}^d)$ у просторі $B_{q,1}(\mathbb{T}^d)$ при  $1 \leqslant q, p < \infty$.

У роботі отримано доповнення і узагальнення деяких результатів одержаних в згаданих вище роботах \cite{Temlyakov-1990Proc-189, Kashin-Temlyakov-1994MN, Belinskii-1998JAT, RomanyukAV-19UMJ2, RomanyukAV-19UMJ8, RomanyukAV-21UMB, Hembarska-Zaderei-2022, Romanyuk-Yanchenko-2021UMJ, Romanyuk-Yanchenko-2022-6UMJ}. При цьому     також  було виявлено специфіку багатовимірної апроксимації класів $W^{\boldsymbol{r}}_{p,\boldsymbol{\alpha}}(\mathbb{T}^d)$  у просторі  $B_{q,1}(\mathbb{T}^d)$, у   порівнянні з   $L_q(\mathbb{T}^d)$-простором. Більш детально про це буде йти мова в коментарях до одержаних  результатів.

\medskip

\textbf{2. Означення функціональних класів}. Нехай $\mathbb{R}^d$, $d\geqslant 1$,~---  евклідів  простір з  елементами $\boldsymbol{x}=(x_1,\ldots, x_d)$ і ${(\boldsymbol{x},\boldsymbol{y}) = x_1 y_1+ \ldots + x_d y_d}$. Через $L_p(\mathbb{T}^d)$, $\mathbb{T}^d = \prod\limits^{d}_{j=1} [0,2\pi)$, $1\leqslant p \leqslant \infty$, позначимо простір функцій $f$, які є  $2\pi$-періодичними за кожною змінною зі скінченною нормою:
  $$
  \|f\|_p : = \|f\|_{L_p(\mathbb{T}^d)} =  \left((2\pi)^{-d}
  \int\limits_{\mathbb{T}^d}|f(\boldsymbol{x})|^p\,d\boldsymbol{x}\right)^{1/p}, \ 1 \leqslant p < \infty,
 $$
 $$
 \|f\|_{\infty}:=\|f\|_{L_{\infty}(\mathbb{T}^d)}=\mathop {\rm ess \sup}\limits_{\boldsymbol{x}\in \mathbb{T}^d} |f(\boldsymbol{x})|.
 $$

Також будемо розглядати множину функцій $L^0_p(\mathbb{T}^d)$, яка визначається таким чином:
$$
L^0_p(\mathbb{T}^d):= \Big\{f\colon f \in L_p(\mathbb{T}^d), \int\limits^{2\pi}_0 f(\boldsymbol{x})dx_j = 0,  j=\overline{1,d}, \ \text{майже скрізь} \Big\}.
$$

Надалі, з метою спрощення записів, також будемо використовувати позначення $L_p$ замість $L_p(\mathbb{T}^d)$ і відповідно $L^0_p$ замість $L^0_p(\mathbb{T}^d)$.

Нехай $F_{\boldsymbol{r}}(\boldsymbol{x}, \boldsymbol{\alpha})$~--- багатовимірні аналоги ядер Бернуллі, тобто
$$
F_{\boldsymbol{r}}({\boldsymbol{x}}, \boldsymbol{\alpha}) = 2^d \sum\limits_{\boldsymbol{k}}\,\prod\limits^d_{j=1} k^{-r_j}_j  \cos\left(k_j x_j - \frac{\alpha_j \pi}{2}\right), r_j>0, \alpha_j \in \mathbb{R},
$$
і підсумовування проводиться за векторами  $\boldsymbol{k}= (k_1, \ldots, k_d)$, для яких $k_j > 0$, $j = \overline{1,d}$. Тоді через $W^{\boldsymbol{r}}_{p,\boldsymbol{\alpha}}(\mathbb{T}^d)$ позначимо клас функцій $f$, які подаються у вигляді
\begin{equation}\label{f-r-alpha-W}
f(\boldsymbol{x}) = \varphi(\boldsymbol{x}) \ast F_{\boldsymbol{r}} (\boldsymbol{x}, \boldsymbol{\alpha}) = (2\pi)^{-d}\int\limits_{\mathbb{T}^d}\varphi(\boldsymbol{y})F_{\boldsymbol{r}}(\boldsymbol{x}-\boldsymbol{y}, \boldsymbol{\alpha})d\boldsymbol{y},
\end{equation}
$$
\varphi \in L_p(\mathbb{T}^d), \|\varphi\|_p \leqslant 1,
$$
де ``$\ast$'' означає операцію згортки.

Функцію $\varphi$ у представленні \eqref{f-r-alpha-W} називають  $(\boldsymbol{r}, \boldsymbol{\alpha})$-похідною  функції $f$ і   позначають $f^{(\boldsymbol{r})}_{\boldsymbol{\alpha}}$.  Зауважимо,    що у випадку $\boldsymbol{\alpha} = \boldsymbol{0}$  класи $W^{\boldsymbol{r}}_{p,\boldsymbol{0}}(\mathbb{T}^d)$  будемо  позначати $W^{\boldsymbol{r}}_{p}(\mathbb{T}^d)$ і відповідно $f^{(\boldsymbol{r})}_{\boldsymbol{0}}: = f^{(\boldsymbol{r})}$.


При доведенні  оцінок  зверху  апроксимаційних характеристик класів $W^{\boldsymbol{r}}_{p,\boldsymbol{\alpha}}(\mathbb{T}^d)$ нам буде зручно отримувати їх для більш широких класів Нікольського--Бєсова $B^{\boldsymbol{r}}_{p,\theta}(\mathbb{T}^d)$ при певних значеннях параметра  $\theta$. Нагадаємо   означення   цих класів у термінах так званого декомпозиційного нормування  (див. \cite{Lizorkin-Nikolsky-1989}  Зауваження 2.1).

Нехай $V_l(t)$,  $t\in\mathbb{R}$, $l \in \mathbb{N}$, позначає ядро Валле Пуссена вигляду
$$
V_l(t) = 1 + 2 \sum\limits^l_{k=1} \cos kt  + 2\sum\limits^{2l-1}_{k=l+1}\left(1 - \frac{k-l}{l}\right)\cos kt,
$$
де  при $l=1$ третій доданок вважаємо рівним нулеві.
Кожному вектору ${\boldsymbol{s} = (s_1, \ldots, s_d)}$, $s_j \in \mathbb{N}$, $j=\overline{1,d}$, поставимо у відповідність поліном
$$
A_{\boldsymbol{s}}(\boldsymbol{x}) = \prod\limits^d_{j=1}\big(V_{2^{s_j}}(x_j) - V_{2^{s_j - 1}}(x_j)\big)
$$
 і для $f \in L^0_p$, $1 \leqslant p \leqslant \infty$,  покладемо
 $$
 A_{\boldsymbol{s}}(f) := A_{\boldsymbol{s}}(f,\boldsymbol{x}) = (f \ast A_{\boldsymbol{s}})(\boldsymbol{x}).
 $$

Тоді при $1 \leqslant p \leqslant \infty$,  $\boldsymbol{r} = (r_1, \ldots, r_d)$, $r_j > 0$, $j = \overline{1,d}$, класи $B^{\boldsymbol{r}}_{p,\theta}(\mathbb{T}^d)$  можна означити таким чином:
$$
B^{\boldsymbol{r}}_{p,\theta}(\mathbb{T}^d):=\Big \{f \in
 L^0_p\colon \|f\|_{B^{\boldsymbol{r}}_{p,\theta}(\mathbb{T}^d)} \leqslant 1 \Big\},
$$
де
$$
\|f\|_{B^{\boldsymbol{r}}_{p,\theta}(\mathbb{T}^d)} \asymp \left(\sum\limits_{\boldsymbol{s}\in \mathbb{N}^d} 2^{(\boldsymbol{s},\boldsymbol{r})\theta} \|A_{\boldsymbol{s}}(f)\|^{\theta}_p \right)^{\frac{1}{\theta}}
$$
при $1 \leqslant \theta < \infty$  і
$$
 \|f\|_{B^{\boldsymbol{r}}_{p,\infty}(\mathbb{T}^d)}\equiv\|f\|_{H^{\boldsymbol{r}}_{p}(\mathbb{T}^d)}\asymp  \sup\limits_{s\in \mathbb{N}^d} 2^{(\boldsymbol{s},\boldsymbol{r})} \|A_{\boldsymbol{s}}(f)\|_p.
$$

Тут і надалі по тексту для додатних величин $a$ і  $b$ вживається запис $a\asymp b$, який означає, що існують такі додатні  сталі $C_1$ та $C_2$,
які не залежать від одного істотного параметра у величинах  $a$ і  $b$, що $C_1 a \leqslant b$ (пишемо $a\ll b$) і  $C_2 a \geqslant b$ (пишемо $a\gg b$).  Всі сталі $C_i$, $i=1,2,\dots$, які зустрічаються у роботі, можуть залежати  лише від тих параметрів, що входять в означення класу, метрики, в якій оцінюється похибка
наближення, та розмірності простору $\mathbb{R}^d$.

Зауважимо, що у випадку $1 < p < \infty$ можна записати означення норми з класів   $B^{\boldsymbol{r}}_{p,\theta}(\mathbb{T}^d)$   у дещо іншій формі, а саме з використанням ``блоків'' ряду Фур'є
функції $f$. Для цього нам знадобляться деякі позначення.

Для векторів
$\boldsymbol{s} = (s_1, \ldots, s_d)$,  $s_j \in \mathbb{N}$,  $\boldsymbol{k}= (k_1, \ldots, k_d)$, $k_j \in \mathbb{Z}$, $j=\overline{1,d}$, покладемо
$$
\rho(\boldsymbol{s}) := \big\{ \boldsymbol{k} = (k_1, \ldots, k_d)\colon 2^{s_j - 1} \leqslant |k_j| < 2^{s_j}, j=\overline{1,d} \big\}
$$
 і для $f \in L^0_p$  позначимо
 $$
 \delta_{\boldsymbol{s}}(f):= \delta_{\boldsymbol{s}}(f, \boldsymbol{x}) = \sum\limits_{\boldsymbol{k} \in \rho(\boldsymbol{s})} \widehat{f}(\boldsymbol{k}) e^{i(\boldsymbol{k},\boldsymbol{x})},
 $$
де
$\widehat{f}(\boldsymbol{k}) =  {\displaystyle\int\limits_{\mathbb{T}^d}} f(\boldsymbol{t}) e^{-i(\boldsymbol{k},\boldsymbol{t})}d\boldsymbol{t}$~--- коефіцієнти Фур'є функції $f$.

Нехай $1<p<\infty$, $1\leqslant  \theta \leqslant \infty$, $\boldsymbol{r}= (r_1, \ldots, r_d)$, $r_j>0$, $j = \overline{1,d}$. Тоді класи $B_{p,\theta}^{\boldsymbol{r}}(\mathbb{T}^d)$ можна означити таким
чином~\cite{Lizorkin-Nikolsky-1989}:
$$
B_{p,\theta}^{\boldsymbol{r}}(\mathbb{T}^d):=\Big\{f\in L^0_p \colon \|f\|_{B^{\boldsymbol{r}}_{p,\theta}(\mathbb{T}^d)} \leqslant 1 \Big\},
$$
де
\begin{equation}\label{Brpt-period-dek1-bs}
   \|f\|_{B^{\boldsymbol{r}}_{p,\theta}(\mathbb{T}^d)}\asymp \left(\sum \limits_{\boldsymbol{s}\in\mathbb{N}^d}2^{(\boldsymbol{s},\boldsymbol{r})\theta}
   \|\delta_{\boldsymbol{s}}(f)\|_p^{\theta}\right)^{\frac{1}{\theta}}
\end{equation}
  при $1\leqslant\theta<\infty$ і
$$
   \|f\|_{B^{\boldsymbol{r}}_{p,\infty}(\mathbb{T}^d)}\equiv\|f\|_{H^{\boldsymbol{r}}_{p}(\mathbb{T}^d)}\asymp \sup
   \limits_{\boldsymbol{s}\in\mathbb{N}^d} 2^{(\boldsymbol{s},\boldsymbol{r})}\|\delta_{\boldsymbol{s}}(f)\|_p.
$$

Нагадаємо, що для введених класів справджуються такі вкладення:
$$
 B^{\boldsymbol{r}}_{p, p}(\mathbb{T}^d) \subset W^{\boldsymbol{r}}_{p, \boldsymbol{\alpha}}(\mathbb{T}^d) \subset B^{\boldsymbol{r}}_{p, 2}(\mathbb{T}^d), \quad 1 < p \leqslant 2;
$$
\begin{equation}\label{Brp-Wr-vklad}
 B^{\boldsymbol{r}}_{p, 2}(\mathbb{T}^d) \subset W^{\boldsymbol{r}}_{p, \boldsymbol{\alpha}}(\mathbb{T}^d) \subset B^{\boldsymbol{r}}_{p, p}(\mathbb{T}^d),  \quad 2 \leqslant p < \infty;
\end{equation}
$$
 W^{\boldsymbol{r}}_{p, \boldsymbol{\alpha}}(\mathbb{T}^d) \subset B^{\boldsymbol{r}}_{p, \infty}(\mathbb{T}^d) \equiv H^{\boldsymbol{r}}_{p}(\mathbb{T}^d), \quad 1 \leqslant p \leqslant \infty.
$$
Зокрема при $ \theta = p = 2$
$$
W^{\boldsymbol{r}}_{2, \boldsymbol{\alpha}}(\mathbb{T}^d) \subset B^{\boldsymbol{r}}_{2, 2}(\mathbb{T}^d) \subset  W^{\boldsymbol{r}}_{2, \boldsymbol{\alpha}}(\mathbb{T}^d).
$$

Надалі замість $W^{\boldsymbol{r}}_{p, \boldsymbol{\alpha}}(\mathbb{T}^d)$, $H^{\boldsymbol{r}}_{p}(\mathbb{T}^d)$ і  $B^{\boldsymbol{r}}_{p, \theta}(\mathbb{T}^d)$  для спрощення записів також будемо писати $W^{\boldsymbol{r}}_{p, \boldsymbol{\alpha}}$, $H^{\boldsymbol{r}}_{p}$ і  $B^{\boldsymbol{r}}_{p, \theta}$, це не повинно створювати непорозумінь, оскільки в роботі розглядаються лише такі класи функцій.

У подальших міркуваннях будемо вважати, що координати векторів
${\boldsymbol{r}=(r_1,\dots,r_d)}$, які входять в означення класів, впорядковані таким чином:
$0<r_1=r_2=\ldots=r_{\nu}<r_{\nu+1}\leqslant\ldots\leqslant r_d$.
Вектору $\boldsymbol{r}=(r_1,\dots,r_d)$ поставимо у відповідність вектор
$\boldsymbol{\gamma}=(\gamma_1,\dots,\gamma_d)$, $\gamma_j={\displaystyle\frac{r_j}{r_1}}$, $j=\overline{1,d}$, а вектору $\boldsymbol{\gamma}$, в свою чергу, вектор $\boldsymbol{\gamma}'=(\gamma_1',\dots,\gamma_d')$, де
$\gamma'_j=\gamma_j$, якщо $j=\overline{1,\nu}$ i
$1<\gamma_j'<\gamma_j$, $j=\overline{\nu+1,d}$.

Тепер наведемо означення норми у підпросторі $B_{q,1}(\mathbb{T}^d)$, $1\leqslant q \leqslant \infty$ (далі пишемо $B_{q,1}$) функцій  $f \in L_q(\mathbb{T}^d)$.

Для тригонометричних  поліномів $t$ за кратною тригонометричною  системою $\{e^{i(\boldsymbol{k},\boldsymbol{x})}\}_{\boldsymbol{k} \in \mathbb{Z}^d}$ норма $\|t\|_{B_{q,1}}$ визначається згідно з формулою
$$
\|t\|_{B_{q,1}} := \sum\limits_{\boldsymbol{s}\in \mathbb{N}^d} \|A_{\boldsymbol{s}}(t)\|_{q}.
$$
Аналогічно означається норма $\|f\|_{B_{q,1}}$  для будь-якої  функції $f \in L_{q}$  такої, що ряд $\sum\limits_{\boldsymbol{s}\in \mathbb{N}^d} \|A_{\boldsymbol{s}}(f)\|_q$  збігається.

Зауважимо, що у випадку $1 < q < \infty$
\begin{equation}\label{Bq1-norm}
\|f\|_{B_{q,1}} \asymp \sum\limits_{\boldsymbol{s} \in \mathbb{N}^d}\|\delta_{\boldsymbol{s}}(f)\|_q.
\end{equation}

Крім того, для $f \in B_{q,1}$, $1 \leqslant q \leqslant \infty$, виконуються співвідношення:
\begin{equation}\label{fq-fBq1}
\|f\|_{q} \ll \|f\|_{B_{q,1}} \quad \text{і}  \quad \|f\|_{B_{1,1}} \ll \|f\|_{B_{q,1}} \ll \|f\|_{B_{\infty, 1}}.
\end{equation}

З історією дослідження різних апроксимаційних характеристик класів $W^{\boldsymbol{r}}_{p,\boldsymbol{\alpha}}(\mathbb{T}^d)$, $H^{\boldsymbol{r}}_p(\mathbb{T}^d)$ і $B^{\boldsymbol{r}}_{p,\theta}(\mathbb{T}^d)$, ${1 \leqslant \theta < \infty}$, у просторах $L_q(\mathbb{T}^d)$, ${1 \leqslant q \leqslant \infty}$, можна ознайомитися у монографіях \cite{Romanyuk-2012m, Temlyakov-1986m, Temlyakov-1993m, Cross-2018, Temlyakov-2018, Trigub-Belinsky} і, зокрема, у просторах $B_{q,1}(\mathbb{T}^d)$, $q \in \{1, \infty\}$~---  у роботах \cite{Temlyakov-1990Proc-189, Kashin-Temlyakov-1994MN, Belinskii-1998JAT, RomanyukAV-19UMJ2, RomanyukAV-19UMJ8, RomanyukAV-21UMB, Hembarska-Zaderei-2022, Romanyuk-Yanchenko-2021UMJ, Romanyuk-Yanchenko-2022-6UMJ, Fedunyk-Hembarska-2020-CMP, RomanyukAV-Pozharska-2023-CMP, Hembarska-Fedunyk-2023}.

\medskip

\textbf{3. Апроксимаційні характеристики і допоміжні твердження}.

Для   $n \in \mathbb{N}$,  $ \boldsymbol{s} \in \mathbb{N}^d$  і $\boldsymbol{\gamma} = (\gamma_1, \ldots, \gamma_d)$, $\gamma_j > 0$, $j = \overline{1,d}$,  позначимо множину
$$
Q^{\boldsymbol{\gamma}}_n := \bigcup\limits_{(\boldsymbol{s},\boldsymbol{\gamma}) < n} \rho(\boldsymbol{s}),
$$
яку називають східчастим  гіперболічним хрестом. У випадку, коли розглядається вектор $\boldsymbol{\gamma'} =(\gamma'_1, \ldots, \gamma'_d)$  відповідну множину  будемо позначати $Q^{\boldsymbol{\gamma'}}_n$ і, зокрема,  при $\boldsymbol{\gamma} = (1, \ldots, 1) \in \mathbb{N}^d$~---  $Q^{\boldsymbol{1}}_n$.

Поставимо  у відповідність множині $Q^{\boldsymbol{\gamma}}_n$ сукупність тригонометричних поліномів вигляду
$$
T(Q^{\boldsymbol{\gamma}}_n):= \Big\{ t\colon t(\boldsymbol{x}) = \sum\limits_{\boldsymbol{k} \in Q^{\boldsymbol{\gamma}}_n} c_{\boldsymbol{k}} e^{i(\boldsymbol{k},\boldsymbol{x})}, \quad c_{\boldsymbol{k}} \in \mathbb{C}, \quad \boldsymbol{x} \in \mathbb{R}^d \Big\}
$$
і для $f \in L^0_1$ покладемо
$$
S_{Q^{\boldsymbol{\gamma}}_n}(f): =S_{Q^{\boldsymbol{\gamma}}_n}(f,\boldsymbol{x}) = \sum\limits_{\boldsymbol{k} \in Q^{\boldsymbol{\gamma}}_n} \widehat{f}(\boldsymbol{k}) e^{i(\boldsymbol{k},\boldsymbol{x})}.
$$

Поліноми $S_{Q^{\boldsymbol{\gamma}}_n}(f)$ називають східчасто-гіперболічними сумами Фур'є функцій  $f$. Згідно з  прийнятими позначеннями, їх також можна подати у вигляді
$$
S_{Q^{\boldsymbol{\gamma}}_n}(f) = \sum\limits_{(\boldsymbol{s}, \boldsymbol{\gamma}) < n} \delta_{\boldsymbol{s}}(f).
$$

Для означених множин тригонометричних поліномів будемо розглядати наступні апроксимаційні характеристики.

Нехай $\mathscr{X}$~--- деякий функціональний простір із нормою $\|\cdot\|_{\mathscr{X}}$. Тоді для $f \in \mathscr{X}$ позначимо через
$$
 E_{Q^{\boldsymbol{\gamma}}_n}\big(f\big)_{\mathscr{X}}:= \inf\limits_{t \in T(Q^{\boldsymbol{\gamma}}_n)} \|f - t\|_{\mathscr{X}}
$$
величину найкращого наближення функції $f$ за допомогою поліномів з множин $T(Q^{\boldsymbol{\gamma}}_n)$.

Відповідно для  функціонального  класу $F\subset \mathscr{X}$   покладемо
\begin{equation}\label{EQngamma}
 E_{Q^{\boldsymbol{\gamma}}_n}\big(F\big)_{\mathscr{X}}:= \sup\limits_{f \in F} E_{Q^{\boldsymbol{\gamma}}_n}\big(f\big)_{\mathscr{X}}.
\end{equation}

Паралельно з величинами \eqref{EQngamma} нами досліджуються також наближення функцій з класів $F$ їхніми східчасто-гіперболічними сумами Фур'є, тобто величини
\begin{equation}\label{EQngamma-Sum}
\mathscr{E}_{Q^{\boldsymbol{\gamma}}_n}\big(F\big)_{\mathscr{X}}:= \sup\limits_{f \in F} \|f - S_{Q^{\boldsymbol{\gamma}}_n}\big(f\big) \|_{\mathscr{X}}.
\end{equation}

Більш конкретно в роботі ми встановлюємо точні за порядком оцінки величин \eqref{EQngamma} і \eqref{EQngamma-Sum} (інколи при $\boldsymbol{\gamma}= \boldsymbol{\gamma'}=(\gamma'_1,\ldots,\gamma'_d)$) для класів
$F = W^{\boldsymbol{r}}_{p,\boldsymbol{\alpha}}$ і просторів $\mathscr{X} = B_{q,1}$, $1 \leqslant p, q < \infty$.

У коментарях до одержаних результатів ми будемо звертатися до одновимірного випадку, тому наведемо відповідні  модифікації означень  величин \eqref{EQngamma} і \eqref{EQngamma-Sum}.

Для $F \subset \mathscr{X}$ позначимо
$$
E_{2^n}\big(F\big)_{\mathscr{X}} : = \sup\limits_{f \in F} \inf\limits_{t \in T(2^n)} \| f - t \|_{\mathscr{X}},
$$
де
$$
T(2^n) := \Big\{ t \colon t(x) = \sum\limits^{2^n}_{k = -2^n} c_k e^{i k x}, \quad c_k \in \mathbb{C}, \quad x \in \mathbb{R} \Big\}.
$$

Відповідно до означення  величини \eqref{EQngamma-Sum} покладемо
$$
\mathscr{E}_{2^n}\big(F\big)_{\mathscr{X}}:= \sup\limits_{f \in F} \|f - S_{2^n}(f)\|_{\mathscr{X}},
$$
де
$$
S_{2^n}(f): = S_{2^n}(f,x):= \sum\limits^{2^n}_{k = -2^n} \widehat{f}(k) e^{i k x}, \quad x\in \mathbb{R}.
$$

Перед тим, як безпосередньо перейти до викладу одержаних результатів, відмітимо важливе співвідношення між   означеними величинами у випадку $f \in  B_{q,1}$, $1 < q< \infty$.

Нехай $\mathbb{S}_{Q^{\boldsymbol{\gamma}}_n}$ позначає оператор Фур'є, який ставить у відповідність функції ${f \in B_{q,1}}$, $1 < q < \infty$, її східчасто-гіперболічну суму Фур'є $S_{Q^{\boldsymbol{\gamma}}_n}(f)$, тобто  $\mathbb{S}_{Q^{\boldsymbol{\gamma}}_n}f = S_{Q^{\boldsymbol{\gamma}}_n}(f)$. Легко  переконатися, що  норма  оператора  $\mathbb{S}_{Q^{\gamma}_n}$ з $B_{q,1}$ в $B_{q,1}$, $1 < q < \infty$, (позначення $\|\mathbb{S}_{Q^{\boldsymbol{\gamma}}_n} \|_{B_{q,1} \rightarrow B_{q,1}}$) є обмеженою.

Згідно з означенням маємо
$$
\|\mathbb{S}_{Q^{\boldsymbol{\gamma}}_n}\|_{B_{q,1} \rightarrow B_{q,1}} = \sup \limits_{\|f\|_{B_{q,1}} \leqslant 1} \|S_{Q^{\boldsymbol{\gamma}}_n}(f)\|_{B_{q,1}} \asymp \sup \limits_{\|f\|_{B_{q,1}} \leqslant 1}  \sum\limits_{\boldsymbol{s} \in \mathbb{N}^d} \|\delta_{\boldsymbol{s}}(S_{Q^{\boldsymbol{\gamma}}_n}(f))\|_q =
$$
\begin{equation}\label{SQng-oper}
= \sup \limits_{\|f\|_{B_{q,1}} \leqslant 1}  \sum\limits_{(\boldsymbol{s}, \boldsymbol{\gamma})< n} \|\delta_{\boldsymbol{s}}(f)\|_q \leqslant \sup \limits_{\|f\|_{B_{q,1}} \leqslant 1} \sum\limits_{\boldsymbol{s} \in \mathbb{N}^d}\|\delta_{\boldsymbol{s}}(f)\|_q \leqslant C_3(q).
\end{equation}

Далі, нехай $t^{\ast} \in T(Q^{\boldsymbol{\gamma}}_n)$~--- поліном найкращого наближення функції $f \in B_{q,1}$.   Тоді, з одного боку,  беручи до уваги,  що  $S_{Q^{\boldsymbol{\gamma}}_n}(t^{\ast}) = t^{\ast}$ і скориставшись \eqref{SQng-oper}, можемо  записати
$$
\mathscr{E}_{Q^{\boldsymbol{\gamma}}_n}\big(f\big)_{B_{q,1}} = \|f - S_{Q^{\boldsymbol{\gamma}}_n}(f) \|_{B_{q,1}} = \|f - t^{\ast} + t^{\ast} - S_{Q^{\boldsymbol{\gamma}}_n}(f)\|_{B_{q,1}} \leqslant
$$
$$
\leqslant \|f - t^{\ast}\|_{B_{q,1}} + \|S_{Q^{\boldsymbol{\gamma}}_n}(f) - t^{\ast}\|_{B_{q,1}} = \|f - t^{\ast}\|_{B_{q,1}}  + \|S_{Q^{\boldsymbol{\gamma}}_n}(f - t^{\ast})\|_{B_{q,1}}  \leqslant
$$
$$
\leqslant E_{Q^{\boldsymbol{\gamma}}_n}\big(f\big)_{B_{q,1}} +  \|\mathbb{S}_{Q^{\boldsymbol{\gamma}}_n}\|_{B_{q,1}\rightarrow B_{q,1}} \|f - t^{\ast}\|_{B_{q,1}} \leqslant
$$
\begin{equation}\label{EQngS-EQng}
\leqslant  E_{Q^{\boldsymbol{\gamma}}_n}\big(f\big)_{B_{q,1}} + C_3(q) E_{Q^{\boldsymbol{\gamma}}_n}\big(f\big)_{B_{q,1}} = C_4(q) E_{Q^{\boldsymbol{\gamma}}_n}\big(f\big)_{B_{q,1}}
\end{equation}

З іншого боку безпосередньо з означень  величин \eqref{EQngamma} і \eqref{EQngamma-Sum} для $f \in B_{q,1}$  маємо
\begin{equation}\label{EQng-EQngS}
E_{Q^{\boldsymbol{\gamma}}_n}\big(f\big)_{B_{q,1}} \leqslant \mathscr{E}_{Q^{\boldsymbol{\gamma}}_n}\big(f\big)_{B_{q,1}}.
\end{equation}

Співставляючи \eqref{EQngS-EQng} і \eqref{EQng-EQngS} отримуємо
\begin{equation}\label{EQng-asymp}
E_{Q^{\boldsymbol{\gamma}}_n}\big(f\big)_{B_{q,1}} \asymp \mathscr{E}_{Q^{\boldsymbol{\gamma}}_n}\big(f\big)_{B_{q,1}}, \quad 1<q<\infty.
\end{equation}
Зрозуміло, що  співвідношення \eqref{EQng-asymp} справедливе і по відношенню до величин $E_{Q^{\boldsymbol{\gamma'}}_n}\big(f\big)_{B_{q,1}}$, $\mathscr{E}_{Q^{\boldsymbol{\gamma'}}_n}\big(f\big)_{B_{q,1}}$  і, зокрема,~--- $E_{Q^{\boldsymbol{1}}_n}\big(f\big)_{B_{q,1}}$  і $\mathscr{E}_{Q^{\boldsymbol{1}}_n}\big(f\big)_{B_{q,1}}$, а  також $E_{2^n}\big(f\big)_{B_{q,1}}$, $\mathscr{E}_{2^n}\big(f\big)_{B_{q,1}}$.

Нагадаємо  означення ще  однієї  апроксимаційної характеристики, яка не досліджується в роботі, але її відомі оцінки будемо використовувати при встановленні оцінок знизу величин $E_{Q^{\boldsymbol{\gamma}}_n}\big(W^{\boldsymbol{r}}_{p,\boldsymbol{\alpha}}\big)_{B_{q,1}}$ і $\mathscr{E}_{Q^{\boldsymbol{\gamma}}_n}\big(W^{\boldsymbol{r}}_{p,\boldsymbol{\alpha}}\big)_{B_{q,1}}$.

Нехай $\mathscr{Y}$~--- нормований простір із нормою $\|\cdot\|_{\mathscr{Y}}$, $\mathfrak{L}_M(\mathscr{Y})$~--- сукупність підпросторів у просторі $\mathscr{Y}$ розмірності, що не перевищує $M$ і $W$~--- центрально-симетрична множина в $\mathscr{Y}$. Величина
$$
d_M\big(W,\mathscr{Y}\big) := \inf\limits_{L_M \in \mathfrak{L}_M(\mathscr{Y})}\sup\limits_{w \in W}\inf\limits_{u \in L_M}\|w - u\|_{\mathscr{Y}}
$$
називається $M$-вимірним  колмогоровським  поперечником множини $W$ у просторі $\mathscr{Y}$ (колмогоровським  поперечником). Поперечник   $d_M\big(W, \mathscr{X}\big)$ ввів у 1936 році А.\,М.~Колмогоров~\cite{Kolmogorof_1936} і він характеризує апроксимаційні можливості $M$-вимірних підпросторів.

Тепер наведемо деякі допоміжні твердження.

\textbf{Лема А} \cite[Вступ]{Temlyakov-1986m}\textbf{.} \it Справедливе співвідношення
$$
 \sum\limits_{(\boldsymbol{s}, \boldsymbol{\gamma}') \geqslant l} 2 ^{-\beta(\boldsymbol{s}, \boldsymbol{\gamma})} \asymp 2^{-\beta l} l^{\nu - 1}, \beta > 0.
$$ \rm

\vskip 2 mm

\textbf{Теорема А} \cite{Romanyuk-UMJ17-10}\textbf{.} \it Нехай $d \geqslant 1$, $1 < p < \infty$, $r_1 > 0$. Тоді при $\boldsymbol{\alpha} \in \mathbb{R}^d$ справедлива оцінка
$$
 d_M\big(W^{\boldsymbol{r}}_{p,\boldsymbol{\alpha}}, B_{1,1}\big) \asymp M^{- r_1} (\log^{\nu -1 }M)^{r_1 + \frac{1}{2}}.
$$ \rm

\vskip 2 mm

\textbf{Теорема Б}\cite{Nikolsky_51}\textbf{.} \it Нехай $\boldsymbol{n} = (n_1, \ldots, n_d)$, $n_j \in \mathbb{N}$, $j = \overline{1, d}$, і
$$
t(\boldsymbol{x}) = \sum\limits_{|k_j| \leqslant  n_j} c_{\boldsymbol{k}} e^{i(\boldsymbol{k},\boldsymbol{x})}.
$$
Тоді  при $1 \leqslant q  < p < \infty$  справедлива   нерівність
\begin{equation}\label{t-neriv}
\|t\|_p  \leqslant 2^d \prod\limits^d_{j=1} n^{\frac{1}q{} - \frac{1}{p}}_j  \|t\|_q.
\end{equation} \rm

\vskip 2 mm

Нерівність \eqref{t-neriv} встановлена С.\,М.~Нікольським і відома як ``нерівність різних метрик''.

\textbf{Лема Б} \cite[гл.1, \S3]{Temlyakov-1986m}\textbf{.} \it Нехай $1 \leqslant p < q < \infty$   і $f \in L^0_p$.  Тоді
$$
\|f\|^q_q \ll \sum\limits_{\boldsymbol{s} \in \mathbb{N}^d}  \|\delta_{\boldsymbol{s}}(f)\|^q_p 2^{\|\boldsymbol{s}\|_1\left(\frac{1}{p} - \frac{1}{q}\right)q},
$$
де $\|\boldsymbol{s}\|_1 = s_1 + \ldots + s_d$. \rm

\vskip 2 mm

Наступне  твердження  відоме  в  математичній літературі, як нерівність Бернштейна для тригонометричних поліномів з множини $T(Q^{\boldsymbol{1}}_n)$.

\textbf{Теорема В} \cite[гл.1, \S3]{Temlyakov-1986m}\textbf{.} \it Нехай $1 < p < \infty$. Тоді при $r_1 \geqslant 0$ справедливе  співвідношення
$$
\sup\limits_{t \in T(Q^{\boldsymbol{1}}_n)} \|t^{(r)}\|_p / \|t\|_p \asymp 2^{n r_1}.
$$
\rm

\medskip

\textbf{3. Наближення  класів   \boldmath{$W^r_{p, \alpha}$} східчасто-гіперболічними сумами Фур'є і їхні найкращі  наближення у просторі  \boldmath{$B_{q,1}$}}.

Справедливе  твердження.

\textbf{Теорема 1.} \it Нехай $d \geqslant 2$, $1 < p  < \infty$, $r_1 > 0$, $\boldsymbol{\alpha} \in \mathbb{R}^d$. Тоді виконуються  співвідношення
\begin{equation}\label{T1-En-W-Bp1}
E_{Q^{\boldsymbol{\gamma'}}_n}\big(W^{\boldsymbol{r}}_{p,\boldsymbol{\alpha}}\big)_{B_{p,1}} \asymp \mathscr{E}_{Q^{\boldsymbol{\gamma'}}_n}\big(W^{\boldsymbol{r}}_{p,\boldsymbol{\alpha}}\big)_{B_{p,1}} \asymp  2^{-nr_1} n^{(\nu  - 1)\xi},
\end{equation}
де $\xi=\max \left\{\frac{1}{2}, \frac{1}{p'}\right\}$, $\frac{1}{p}+\frac{1}{p'}=1$.
\rm

\emph{\textbf{Доведення.}} Одержимо спочатку оцінку зверху величини $\mathscr{E}_{Q^{\boldsymbol{\gamma'}}_n}\big(W^{\boldsymbol{r}}_{p,\boldsymbol{\alpha}}\big)_{B_{p,1}}$, з  якої, згідно з \eqref{EQng-EQngS} буде випливати оцінка зверху і для найкращого  наближення  $E_{Q^{\boldsymbol{\gamma'}}_n}\big(W^{\boldsymbol{r}}_{p,\boldsymbol{\alpha}}\big)_{B_{p,1}}$. Для цього розглянемо два випадки.

а) Нехай спочатку $p\in (1,2]$. У цьому випадку нам буде зручно одержати шукану оцінку в більш загальній ситуації,  а саме для класів $B^{\boldsymbol{r}}_{p,2}$.

Отже, для довільної функції $f\in B^{\boldsymbol{r}}_{p,2}$ згідно з означенням норми у просторі $B_{p,1}$ \eqref{Bq1-norm} можемо записати
$$
\mathscr{E}_{Q^{\boldsymbol{\gamma'}}_n}\big(f\big)_{B_{p,1}} = \Bigg\|f - \sum\limits_{(\boldsymbol{s},\boldsymbol{\gamma'})<n} \delta_{\boldsymbol{s}}(f)\Bigg\|_{B_{p,1}} = \Bigg\|\sum\limits_{(\boldsymbol{s},\boldsymbol{\gamma'})\geqslant n} \delta_{\boldsymbol{s}}(f)\Bigg\|_{B_{p,1}} \asymp
$$
\begin{equation}\label{EQn-f-T1-1}
\asymp \sum\limits_{\boldsymbol{s} \in \mathbb{N}^d} \left\|\delta_{\boldsymbol{s}}\left(\mathop {\sum_{\boldsymbol{s'} \in \mathbb{N}^d}} \limits_{(\boldsymbol{s'}, \boldsymbol{\gamma'}) \geqslant n} \delta_{\boldsymbol{s'}}(f) \right)\right\|_p \leqslant \sum\limits_{(\boldsymbol{s},\boldsymbol{\gamma'})\geqslant n} \|\delta_{\boldsymbol{s}}(f)\|_p
 = J_1.
\end{equation}

Далі, скориставшись нерівністю Коші--Буняковського (Гельдера з показником $2$), Лемою~А і означенням норми \eqref{Brpt-period-dek1-bs}, отримаємо
$$
 J_1 = \sum\limits_{(\boldsymbol{s},\boldsymbol{\gamma'})\geqslant n} 2^{(\boldsymbol{s},\boldsymbol{r})} \|\delta_{\boldsymbol{s}}(f)\|_p 2^{-(\boldsymbol{s},\boldsymbol{r})} \leqslant  \left(\sum\limits_{(\boldsymbol{s},\boldsymbol{\gamma'})\geqslant n} 2^{2(\boldsymbol{s},\boldsymbol{r})} \|\delta_{\boldsymbol{s}}(f)\|^2_p \right)^{\frac{1}{2}} \times
$$
\begin{equation}\label{T1-J1-os}
\times  \left(\sum\limits_{(\boldsymbol{s},\boldsymbol{\gamma'})\geqslant n} 2^{-2r_1(\boldsymbol{s},\boldsymbol{\gamma})}\right)^{\frac{1}{2}} \ll    \|f\|_{B^{\boldsymbol{r}}_{p,2}} \left(\sum\limits_{(\boldsymbol{s},\boldsymbol{\gamma'})\geqslant n} 2^{-2r_1(\boldsymbol{s},\boldsymbol{\gamma})} \right)^{\frac{1}{2}}  \ll 2^{-n r_1}n^{\frac{\nu - 1}{2}}.
\end{equation}

Тепер врахувавши, що виконується співвідношення \eqref{Brp-Wr-vklad}, з якого при $p \in (1,2]$ маємо  $W^{\boldsymbol{r}}_{p,\boldsymbol{\alpha}} \subset B^{\boldsymbol{r}}_{p,2}$, із співвідношень \eqref{EQn-f-T1-1} і \eqref{T1-J1-os}   знаходимо
\begin{equation}\label{T1-os1}
E_{Q^{\boldsymbol{\gamma'}}_n}\big(W^{\boldsymbol{r}}_{p,\boldsymbol{\alpha}}\big)_{B_{p,1}} \leqslant \mathscr{E}_{Q^{\boldsymbol{\gamma'}}_n}\big(W^{\boldsymbol{r}}_{p,\boldsymbol{\alpha}}\big)_{B_{p,1}} \ll \mathscr{E}_{Q^{\boldsymbol{\gamma'}}_n}\big(B^{\boldsymbol{r}}_{p,2}\big)_{B_{p,1}} \ll 2^{-nr_1} n^{\frac{\nu -1}{2}}.
\end{equation}

б) Нехай $p\in(2;\infty)$. Як і у попередньому випадку, тут нам також буде зручно одержати шукану оцінку в більш загальній ситуації, а саме для класів $B^{\boldsymbol{r}}_{p,p}$. Отже, для $f\in B^{\boldsymbol{r}}_{p,p}$ згідно з означенням норми у просторі $B_{p,1}$, нерівністю Гельдера, означенням норми \eqref{Brpt-period-dek1-bs} і Лемою~А  будемо мати
$$
\mathscr{E}_{Q^{\boldsymbol{\gamma'}}_n}\big(f\big)_{B_{p,1}} = \Bigg\|f - \sum\limits_{(\boldsymbol{s},\boldsymbol{\gamma'})<n} \delta_{\boldsymbol{s}}(f)\Bigg\|_{B_{p,1}} = \Bigg\|\sum\limits_{(\boldsymbol{s},\boldsymbol{\gamma'})\geqslant n} \delta_{\boldsymbol{s}}(f)\Bigg\|_{B_{p,1}} \asymp
$$
$$
\asymp \sum\limits_{\boldsymbol{s} \in \mathbb{N}^d} \left\|\delta_{\boldsymbol{s}}\left(\mathop {\sum_{\boldsymbol{s'} \in \mathbb{N}^d}} \limits_{(\boldsymbol{s'}, \boldsymbol{\gamma'}) \geqslant n} \delta_{\boldsymbol{s'}}(f) \right)\right\|_p \leqslant \sum\limits_{(\boldsymbol{s},\boldsymbol{\gamma'})\geqslant n} \|\delta_{\boldsymbol{s}}(f)\|_p
\leqslant
$$
$$
\leqslant  \left(\sum\limits_{(\boldsymbol{s},\boldsymbol{\gamma'})\geqslant n} 2^{(\boldsymbol{s},\boldsymbol{r})p} \|\delta_{\boldsymbol{s}}(f)\|^p_p \right)^{\frac{1}{p}} \left(\sum\limits_{(\boldsymbol{s},\boldsymbol{\gamma'})\geqslant n} 2^{-(\boldsymbol{s},\boldsymbol{r})p'}\right)^{\frac{1}{p'}} \ll
$$
\begin{equation}\label{T1-2-os}
\ll    \|f\|_{B^{\boldsymbol{r}}_{p,p}} \left(\sum\limits_{(\boldsymbol{s},\boldsymbol{\gamma'})\geqslant n} 2^{-pr_1(\boldsymbol{s},\boldsymbol{\gamma})} \right)^{\frac{1}{p'}}  \ll 2^{-n r_1}n^{(\nu - 1)\frac{1}{p'}}.
\end{equation}

Отже, врахувавши, що $W^{\boldsymbol{r}}_{p,\boldsymbol{\alpha}} \subset B^{\boldsymbol{r}}_{p,p}$, $p \in (2,\infty)$, із \eqref{T1-2-os} одержимо
\begin{equation}\label{T1-os2}
E_{Q^{\boldsymbol{\gamma'}}_n}\big(W^{\boldsymbol{r}}_{p,\boldsymbol{\alpha}}\big)_{B_{p,1}} \leqslant \mathscr{E}_{Q^{\boldsymbol{\gamma'}}_n}\big(W^{\boldsymbol{r}}_{p,\boldsymbol{\alpha}}\big)_{B_{p,1}} \ll \mathscr{E}_{Q^{\boldsymbol{\gamma'}}_n}\big(B^{\boldsymbol{r}}_{p,p}\big)_{B_{p,1}} \ll 2^{-n r_1} n^{(\nu-1)\frac{1}{p'}}.
\end{equation}

Таким чином поєднавши \eqref{T1-os1} і \eqref{T1-os2} отримаємо шукану оцінку зверху.

При встановленні в \eqref{T1-En-W-Bp1} відповідних оцінок знизу також будемо розглядати два випадки.

а) Нехай $p\in (1,2]$. Тоді оцінки знизу для обох апроксимаційних характеристик є наслідками відомої оцінки колмогоровського поперечника $d_M\big(W^{\boldsymbol{r}}_{p,\boldsymbol{\alpha}}, B_{1,1}\big)$ (Теорема~А).

Так, вибравши число $n\in \mathbb{N}$ із співвідношення $M\asymp 2^n n^{\nu-1}$ отримаємо
$$
\mathscr{E}_{Q^{\boldsymbol{\gamma'}}_n}\big(W^{\boldsymbol{r}}_{p,\boldsymbol{\alpha}}\big)_{B_{p,1}} \geqslant E_{Q^{\boldsymbol{\gamma'}}_n}\big(W^{\boldsymbol{r}}_{p,\boldsymbol{\alpha}}\big)_{B_{p,1}} \geqslant E_{Q^{\boldsymbol{\gamma'}}_n}\big(W^{\boldsymbol{r}}_{p,\boldsymbol{\alpha}}\big)_{B_{1,1}} \gg
$$
$$
\gg d_M\big(W^{\boldsymbol{r}}_{p,\boldsymbol{\alpha}}, B_{1,1}\big) \asymp  M^{-r_1}\big(\log^{\nu - 1}M\big)^{r_1+\frac{1}{2}} \asymp   2^{-nr_1} n^{\frac{\nu  - 1}{2}}.
$$

б) Нехай $p\in (2,\infty)$. Зауважимо, що необхідні оцінки знизу достатньо встановити при $\boldsymbol{\alpha}=\boldsymbol{0}$ і $\nu=d$. Розглянемо функцію $g_1 \in W^{\boldsymbol{r}}_p$, $\boldsymbol{r} = (r_1, \ldots, r_1) \in \mathbb{R}^d$, наближення якої у просторі $B_{p,1}$ східчасто-гіперболічними сумами Фур'є $S_{Q^{\boldsymbol{1}}_n}(g_1)$ буде реалізувати одержані вище порядкові оцінки зверху.

Нехай
$$
g_1(\boldsymbol{x}) = C_5 2^{-n\left(r_1 + 1 - \frac{1}{p}\right)} n^{-\frac{d - 1}{p}} d_n(\boldsymbol{x}), \quad C_5 >0,
$$
де
$$
d_n(\boldsymbol{x}) = \sum\limits_{(\boldsymbol{s},\boldsymbol{1}) = n} \sum\limits_{\boldsymbol{k} \in \rho(\boldsymbol{s})} e^{i(\boldsymbol{k},\boldsymbol{x})}
$$
і покажемо, що при відповідному виборі сталої $C_5 > 0$ функція  $g_1$ належить класу $W^{\boldsymbol{r}}_p$.

Для цього скористаємося спочатку лемою~Б, згідно з якою при $p \in (2,\infty)$ можемо записати
$$
\|g_1\|_p \ll  \left(\sum\limits_{(\boldsymbol{s},\boldsymbol{1}) = n} \|\delta_{\boldsymbol{s}}(g_1)\|^p_2 2^{\|\boldsymbol{s}\|_1 \left(\frac{1}{2} - \frac{1}{p}\right)p}\right)^{\frac{1}{p}} \asymp
$$
$$
\asymp 2^{-n\left(r_1+1-\frac{1}{p}\right)} n^{-\frac{d - 1}{p}} \left(\sum\limits_{(\boldsymbol{s},\boldsymbol{1}) = n} \|\delta_{\boldsymbol{s}}(d_n)\|^p_2  2^{\|\boldsymbol{s}\|_1\left(\frac{1}{2} - \frac{1}{p}\right)p}   \right)^{\frac{1}{p}} = J_2.
$$

Далі скориставшись співвідношенням (див., наприклад, \cite[гл.~1, \S1]{Temlyakov-1993m})
$$
\Bigg\|\sum\limits^m_{k = -m} e^{i k x}\Bigg\|_q \asymp m^{1 - \frac{1}{q}}, \quad q \in (1, \infty),
$$
маємо
\begin{equation}\label{dn-zag}
\Bigg\|\sum\limits_{\boldsymbol{k} \in \rho(\boldsymbol{s})} e^{i(\boldsymbol{k},\boldsymbol{x})}\Bigg\|_q \asymp 2^{\|\boldsymbol{s}\|_1\left(1 - \frac{1}{q}\right)}, \quad q \in (1,\infty).
\end{equation}
Використавши \eqref{dn-zag} при $q=2$ продовжимо оцінку величини $J_2$:
$$
J_2 \asymp 2^{-n\left(r_1+1-\frac{1}{p}\right)} n^{-\frac{d-1}{p}} \left(\sum\limits_{(\boldsymbol{s},\boldsymbol{1})=n} 2^{\|\boldsymbol{s}\|_1 \left(1-\frac{1}{p}\right)p}\right)^{\frac{1}{p}} \asymp
$$
\begin{equation}\label{T1-J2-os-N2}
\asymp 2^{-n\left(r_1+1-\frac{1}{p}\right)} n^{-\frac{d - 1}{p}} 2^{n\left(1-\frac{1}{p}\right)} n^{\frac{d - 1}{p}} = 2^{-n r_1}.
\end{equation}

На завершення, скориставшись теоремою~В, із \eqref{T1-J2-os-N2} отримуємо
\begin{equation}\label{T1-g1-os}
\|g^{(\boldsymbol{r})}_1\|_p \ll  2^{n r_1} \|g_1\|_p \ll 2^{n r_1} 2^{-n r_1} = 1.
\end{equation}
Останнє співвідношення показує, що   $g_1 \in W^{\boldsymbol{r}}_p$, $p \in (2, \infty)$.

Зазначимо, що співвідношення аналогічне до \eqref{T1-g1-os} справедливе і у випадку $p=2$.

Для продовження міркувань зауважимо, що  $S_{Q^{\boldsymbol{1}}_n}(g_1) = 0$, і тому будемо мати
$$
E_{Q^{\boldsymbol{1}}_n}\big(g_1\big)_{B_{p,1}} \asymp \mathscr{E}_{Q^{\boldsymbol{1}}_n}\big(g_1\big)_{B_{p,1}} = \|g_1\|_{B_{p,1}}\asymp 2^{-n \left(r_1+1-\frac{1}{p}\right)} n^{-\frac{d-1}{p}} \sum\limits_{(\boldsymbol{s},\boldsymbol{1})=n} \|\delta_{\boldsymbol{s}}(d_n)\|_p \asymp
$$
$$
\asymp 2^{-n\left(r_1+1-\frac{1}{p}\right)} n^{-\frac{d - 1}{p}} \sum\limits_{(\boldsymbol{s},\boldsymbol{1})=n} 2^{(\boldsymbol{s},\boldsymbol{1})\left(1-\frac{1}{p}\right)} \asymp
$$
$$
\asymp 2^{-n\left(r_1+1- \frac{1}{p}\right)} n^{-\frac{d - 1}{p}} 2^{n\left(1-\frac{1}{p}\right)} n^{d-1} =  2^{-nr_1} n^{(d - 1)\left(1 - \frac{1}{p}\right)}=  2^{-nr_1} n^{(d - 1)\frac{1}{p'}}.
$$

Оцінки знизу встановлено.

Теорему~1 доведено.

\vskip 1 mm 

\emph{\textbf{Зауваження 1.}} За допомогою аналогічних міркувань легко переконатися, що при виконанні умов теореми~1 виконуються співвідношення
\begin{equation}\label{T1-En-W-Bp1-Zayv}
E_{Q^{\boldsymbol{\gamma}}_n}\big(W^{\boldsymbol{r}}_{p,\boldsymbol{\alpha}}\big)_{B_{p,1}} \asymp \mathscr{E}_{Q^{\boldsymbol{\gamma}}_n}\big(W^{\boldsymbol{r}}_{p,\boldsymbol{\alpha}}\big)_{B_{p,1}} \asymp  2^{-nr_1} n^{(d  - 1)\xi},
\end{equation}
\vskip 2mm

В  доповнення   до   теореми 1 наведемо   твердження стосовно   одновимірного випадку.

\textbf{Теорема~1$'$.} \it Нехай $d =1$, $1<p<\infty$, $r>0$, $\alpha \in \mathbb{R}$. Тоді   справедливі  оцінки
\begin{equation}\label{T1'-En-W-Bp1-d=1}
E_{2^n}\big(W^{r}_{p,\alpha}\big)_{B_{p,1}} \asymp \mathscr{E}_{2^n}\big(W^{r}_{p,\alpha}\big)_{B_{p,1}} \asymp  2^{-n r}.
\end{equation} \rm

\emph{\textbf{Доведення.}} Оцінки зверху в \eqref{T1'-En-W-Bp1-d=1} випливають з наслідку~1 \cite{RomanyukAV-Pozharska-2023-CMP} при $d = 1$ згідно зі вкладенням $W^r_{p,\alpha} \subset H^{r}_p$ і співвідношеннями
$$
E_{2^n}\big(W^{r}_{p,\alpha}\big)_{B_{p,1}} \leqslant \mathscr{E}_{2^n}\big(W^{r}_{p,\alpha}\big)_{B_{p,1}} \ll \mathscr{E}_{2^n}\big(H^r_p\big)_{B_{p,1}} \asymp 2^{-n r}.
$$

Оцінки знизу в \eqref{T1'-En-W-Bp1-d=1} є наслідком оцінок відповідних величин у просторі $L_p$ (див., наприклад, \cite[гл.~1, \S\,3]{Temlyakov-1993m} і порядкової нерівності \eqref{fq-fBq1} ($\|\cdot\|_{B_{p,1}} \gg \|\cdot\|_p$).

Теорему~1$'$  доведено. \vskip 2mm

Як наслідок отриманого результату, можемо зробити  висновок, що в одновимірному  випадку відповідні  апроксимаційні характеристики класів $W^{r}_{p,\alpha}$ у просторах $B_{p,1}$ і $L_p$ співпадають за порядком.

У випадку $d \geqslant 2$ ситуація є іншою. Наведемо твердження, яке відповідає теоремі~1 у просторі $L_p$.

\textbf{Теорема Г.} \it Нехай $d \geqslant 2$, $1<p<\infty$, $r_1>0$. Тоді при $\boldsymbol{\alpha} \in \mathbb{R}^d$ справедливі  оцінки
\begin{equation}\label{T-En-W-Lp}
E_{Q^{\boldsymbol{\gamma}}_n}\big(W^{\boldsymbol{r}}_{p,\boldsymbol{\alpha}}\big)_p \asymp   \mathscr{E}_{Q^{\boldsymbol{\gamma}}_n}\big(W^{\boldsymbol{r}}_{p,\boldsymbol{\alpha}}\big)_p \asymp 2^{-n r_1}.
\end{equation}\rm \vskip 2 mm

Нагадаємо, що у  випадку, коли $\boldsymbol{r}$~--- вектор з цілочисельними координатами оцінки встановлені в \cite{Mityagin-1962}, а для довільного вектора $\boldsymbol{r}$~--- в \cite{Nikolskaya-75-sib}.

У зв'язку з оцінкою \eqref{T-En-W-Lp} відмітимо також справедливість співвідношень
\begin{equation}\label{T-En'-W-Lp}
E_{Q^{\boldsymbol{\gamma'}}_n}\big(W^{\boldsymbol{r}}_{p,\boldsymbol{\alpha}}\big)_p \asymp   \mathscr{E}_{Q^{\boldsymbol{\gamma'}}_n}\big(W^{\boldsymbol{r}}_{p,\boldsymbol{\alpha}})_p \asymp 2^{-n r_1}.
\end{equation}

Таким чином, співставивши \eqref{T1-En-W-Bp1} і \eqref{T1-En-W-Bp1-Zayv}  з \eqref{T-En'-W-Lp} і \eqref{T-En-W-Lp} виявляємо, що при  $d \geqslant 2$ порядки відповідних апроксимаційних характеристик класів $W^{\boldsymbol{r}}_{p,\boldsymbol{\alpha}}$ є однаковими лише при $\nu = 1$. Крім цього у просторі $B_{p,1}$ існує залежність оцінок цих характеристик як від показників розмірності $\nu$ або $d$, так і від значення параметра $p$. Більше того, при $\nu \neq d$ порядки наближень поліномами з множин $T(Q^{\boldsymbol{\gamma}}_n)$ і $T(Q^{\boldsymbol{\gamma'}}_n)$ у просторі $B_{p,1}$ є різними на відміну від  наближень в $L_p$-просторі.

У наступному твердженні доповнимо оцінки \eqref{T1-En-W-Bp1}, розглянувши випадок $p=1$, але лише для найкращих наближень $E_{Q^{\boldsymbol{\gamma'}}_n}\big(W^{\boldsymbol{r}}_{1,\boldsymbol{\alpha}}\big)_{B_{1,1}}$.

\textbf{Теорема 2.} \it Нехай $d \geqslant 2$, $r_1 > 0$  і $\boldsymbol{\alpha} \in \mathbb{R}^d$.  Тоді  справедлива оцінка
\begin{equation}\label{T2-En-W-B11}
 E_{Q^{\boldsymbol{\gamma'}}_n}\big(W^{\boldsymbol{r}}_{1,\boldsymbol{\alpha}}\big)_{B_{1,1}} \asymp  2^{ -n  r_1} n^{\nu - 1}.
\end{equation}\rm

\emph{\textbf{Доведення.}} Встановимо в \eqref{T2-En-W-B11} оцінку  зверху для класів $H^{\boldsymbol{r}}_1$, з якої згідно із вкладенням $W^{\boldsymbol{r}}_{1,\boldsymbol{\alpha}} \subset H^{\boldsymbol{r}}_1$ буде випливати шукана оцінка і для класів $W^{\boldsymbol{r}}_{1,\boldsymbol{\alpha}}$.

Нехай $f \in H^{\boldsymbol{r}}_1$. В якості агрегату наближення функції $f$ розглянемо поліном вигляду
$$
t_n:= t_n(\boldsymbol{x})= \sum\limits_{(\boldsymbol{s},\boldsymbol{\gamma'})< n-\boldsymbol{\gamma'}(d)} A_{\boldsymbol{s}}(f,\boldsymbol{x}),
$$
де $\boldsymbol{\gamma'}(d) = \gamma'_1+\ldots+\gamma'_d$, $n>3\boldsymbol{\gamma'}(d)$. Легко бачити, що  $t_n \in T(Q^{\boldsymbol{\gamma'}}_n)$.

Тоді на підставі означення норми у просторі $B_{1,1}$ можемо записати
$$
E_{Q^{\boldsymbol{\gamma'}}_n}\big(f\big)_{B_{1,1}} \leqslant  \|f-t_n\|_{B_{1,1}} =\Bigg\|\sum\limits_{(\boldsymbol{s}, \boldsymbol{\gamma'}) \geqslant n-\boldsymbol{\gamma'}(d)} A_{\boldsymbol{s}}(f)\Bigg\|_{B_{1,1}} =
$$
$$
= \sum\limits_{\boldsymbol{s} \in \mathbb{N}^d} \Bigg\|A_{\boldsymbol{s}} \ast \mathop{\sum_{\boldsymbol{s'} \in \mathbb{N}^d}} \limits_{(\boldsymbol{s'},\boldsymbol{\gamma'}) \geqslant n-\boldsymbol{\gamma'}(d)}  A_{\boldsymbol{s'}}(f)\Bigg\|_1 \leqslant \sum\limits_{(\boldsymbol{s},\boldsymbol{\gamma'}) \geqslant n-2\boldsymbol{\gamma'}(d)} \Bigg\|A_{\boldsymbol{s}} \ast \sum\limits_{\|\boldsymbol{s}-\boldsymbol{s'}\|_{\infty}\leqslant 1} A_{\boldsymbol{s'}}(f)\Bigg\|_1 \leqslant
$$
\begin{equation}\label{T2-J3-os1}
\leqslant \sum\limits_{(\boldsymbol{s},\boldsymbol{\gamma'}) \geqslant n-2\boldsymbol{\gamma'}(d)} \|A_{\boldsymbol{s}}\|_1 \Bigg\|\sum\limits_{\|\boldsymbol{s}-\boldsymbol{s'}\|_{\infty} \leqslant 1} A_{\boldsymbol{s'}}(f)\Bigg\|_1 = J_3.
\end{equation}

Для продовження оцінки величини $J_3$ скористаємося тим, що згідно зі співвідношенням $\|V_{2^s}\|_1 \leqslant C_6$, $C_6>0$ (див., наприклад, \cite[гл.\,1, \S\,1]{Temlyakov-1993m})
\begin{equation}\label{As-norn-d1}
\|A_s\|_1 = \|V_{2^s}-V_{2^{s-1}}\|_1 \leqslant \|V_{2^s}\|_1 + \|V_{2^{s-1}}\|_1 \leqslant C_7, \quad C_7 >0.
\end{equation}
Крім цього зауважимо, що для $f \in H^{\boldsymbol{r}}_1$ виконується оцінка
\begin{equation}\label{As-norn-H-d}
 \|A_{\boldsymbol{s}}(f)\|_1 \ll 2^{-(\boldsymbol{s},\boldsymbol{r})}, \quad \boldsymbol{s}\in \mathbb{N}^d.
\end{equation}

Таким чином, із врахуванням \eqref{As-norn-d1}, \eqref{As-norn-H-d} і леми~А, для оцінки величини $J_3$ записуємо
$$
J_3 \ll \sum\limits_{(\boldsymbol{s},\boldsymbol{\gamma'}) \geqslant n-2\boldsymbol{\gamma'}(d)} \sum\limits_{\|\boldsymbol{s}-\boldsymbol{s'}\|_{\infty} \leqslant 1} \|A_{\boldsymbol{s'}}(f)\|_1 \ll
$$
$$
\ll \sum\limits_{(\boldsymbol{s},\boldsymbol{\gamma'}) \geqslant n-3\boldsymbol{\gamma'}(d)} \|A_{\boldsymbol{s}}(f)\|_1 \ll \sum\limits_{(\boldsymbol{s},\boldsymbol{\gamma'}) \geqslant n-3\boldsymbol{\gamma'}(d)} 2^{-(\boldsymbol{s},\boldsymbol{r})} =
$$
\begin{equation}\label{T2-J3-os2}
  = \sum\limits_{(\boldsymbol{s},\boldsymbol{\gamma'}) \geqslant n-3\boldsymbol{\gamma'}(d)} 2^{-(\boldsymbol{s},\boldsymbol{\gamma})r_1} \asymp 2^{-n r_1} n^{\nu-1}.
\end{equation}

Отже, поєднавши  \eqref{T2-J3-os1} і \eqref{T2-J3-os2} та врахувавши, що  $W^{\boldsymbol{r}}_{1,\boldsymbol{\alpha}} \subset H^{\boldsymbol{r}}_1$, отримуємо
$$
E_{Q^{\boldsymbol{\gamma'}}_n}\big(W^{\boldsymbol{r}}_{1,\boldsymbol{\alpha}}\big)_{B_{1,1}} \ll E_{Q^{\boldsymbol{\gamma'}}_n}\big(H^{\boldsymbol{r}}_1\big)_{B_{1,1}} \ll 2^{-n r_1} n^{\nu-1}.
$$

Оцінка зверху встановлена.

Відповідна оцінка знизу в \eqref{T2-En-W-B11} є наслідком співвідношення
$$
E_{Q^{\boldsymbol{\gamma'}}_n}\big(W^{\boldsymbol{r}}_{1,\boldsymbol{\alpha}}\big)_1 \asymp 2^{-n r_1} n^{\nu-1}
$$
(див. \cite[Теорема~4.3]{Temlyakov-1986m}) і нерівності $\|\cdot\|_{B_{1,1}} \gg \|\cdot\|_1$.

Теорему~2 доведено.

\vskip 1 mm 

\emph{\textbf{Зауваження 2.}} За допомогою аналогічних міркувань можна переконатися, що при виконанні умов теореми~2 справедливою є оцінка
$$
E_{Q^{\boldsymbol{\gamma}}_n}\big(W^{\boldsymbol{r}}_{1,\boldsymbol{\alpha}})_{B_{1,1}} \asymp 2^{-n r_1} n^{d-1}.
$$ \vskip 2 mm

Відповідне теоремі~2 твердження в одновимірному випадку має наступний вигляд.

\textbf{Теорема~2$'$.} \it Нехай $d =1$, $r>0$ і $\alpha \in \mathbb{R}$. Тоді
\begin{equation}\label{T2'-En-W-B11-d1}
E_{2^n}\big(W^{r}_{1,\alpha}\big)_{B_{1,1}} \asymp 2^{-n r}.
\end{equation} \rm

\emph{\textbf{Доведення.}} Оцінка зверху випливає з теореми~2 роботи \cite{RomanyukAV-Pozharska-2023-CMP} згідно з вкладенням  $W^r_{1,\alpha} \subset H^r_1$ і співвідношеннями
$$
E_{2^n}\big(W^{r}_{1,\alpha}\big)_{B_{1,1}} \ll E_{2^n}\big(H^r_1\big)_{B_{1,1}} \asymp 2^{-n r}.
$$

Оцінка знизу в \eqref{T2'-En-W-B11-d1} є наслідком оцінки величини $E_{2^n}\big(W^{r}_{1,\alpha}\big)_1$ (див., наприклад, \cite[гл.\,1, \S\,3]{Temlyakov-1993m} і нерівності $\|\cdot\|_{B_{1,1}} \gg \|\cdot\|_1$.

Теорему~2$'$ доведено. \vskip 2 mm

Проаналізувавши результати теорем~2, 2$'$ та їхні доведення можна зробити висновок, що для всіх розмірностей $d\geqslant1$ оцінки найкращих наближень класів $W^{\boldsymbol{r}}_{1,\boldsymbol{\alpha}}$ у просторах $B_{1,1}$ і $L_1$ є однаковими за порядком.

Далі встановимо оцінки величин \eqref{EQngamma} і \eqref{EQngamma-Sum} у випадках, коли відповідні параметри $p$ і $q$ є різними у досліджуваних класах та просторах, в метриці яких оцінюється похибка наближення.

\textbf{Теорема~3.} \it Нехай $d \geqslant 2$, $2 \leqslant p < q < \infty$, $r_1 >\frac{1}{p}-\frac{1}{q}$. Тоді при $\boldsymbol{\alpha}\in \mathbb{R}^d$ справедливі співвідношення
\begin{equation}\label{T3-En-Wp-Bq1}
E_{Q^{\boldsymbol{\gamma}}_n}\big(W^{\boldsymbol{r}}_{p,\boldsymbol{\alpha}}\big)_{B_{q,1}} \asymp   \mathscr{E}_{Q^{\boldsymbol{\gamma}}_n}\big(W^{\boldsymbol{r}}_{p,\boldsymbol{\alpha}}\big)_{B_{q,1}}
 \asymp 2^{-n\left(r_1-\frac{1}{p}+\frac{1}{q}\right)} n^{(\nu-1)\left(1 - \frac{1}{p}\right)}.
\end{equation} \rm

\emph{\textbf{Доведення.}} Оцінку зверху величини $\mathscr{E}_{Q^{\boldsymbol{\gamma}}_n}\big(W^{\boldsymbol{r}}_{p,\boldsymbol{\alpha}}\big)_{B_{q,1}}$, аналогічно до того, як і при доведенні теореми~1, встановимо в більш загальній ситуації, а саме для класів $B^{\boldsymbol{r}}_{p,p}$, оскільки, як зазначалося вище $W^{\boldsymbol{r}}_{p,\boldsymbol{\alpha}} \subset B^{\boldsymbol{r}}_{p,p}$, $2 \leqslant p < \infty$.

Отже, нехай $f$~--- довільна функція із  класу $B^{\boldsymbol{r}}_{p,p}$. Згідно із означенням норми у просторі $B_{q,1}$ і нерівністю \eqref{t-neriv} можемо записати
$$
E_{Q^{\boldsymbol{\gamma}}_n}\big(f\big)_{B_{q,1}}= \Bigg\|f-\sum\limits_{(\boldsymbol{s},\boldsymbol{\gamma})< n} \delta_{\boldsymbol{s}}(f)\bigg\|_{B_{q,1}} = \Bigg\|\sum\limits_{(\boldsymbol{s},\boldsymbol{\gamma}) \geqslant  n} \delta_{\boldsymbol{s}}(f)\Bigg\|_{B_{q,1}} \asymp
$$
$$
\asymp \sum\limits_{\boldsymbol{s} \in \mathbb{N}^d} \left\|\delta_{\boldsymbol{s}} \left(\mathop{\sum_{\boldsymbol{s'} \in \mathbb{N}^d}} \limits_{(\boldsymbol{s'},\boldsymbol{\gamma}) \geqslant n} \delta_{\boldsymbol{s'}}(f)\right)\right\|_q \leqslant \sum\limits_{(\boldsymbol{s},\boldsymbol{\gamma}) \geqslant n} \|\delta_{\boldsymbol{s}}(f)\|_q \ll
$$
\begin{equation}\label{T3-J4-os1}
\ll \sum\limits_{(\boldsymbol{s},\boldsymbol{\gamma}) \geqslant n} 2^{\|\boldsymbol{s}\|_1\left(\frac{1}{p}-\frac{1}{q}\right)} \|\delta_{\boldsymbol{s}}(f)\|_p = J_4.
\end{equation}

Далі, скориставшись нерівністю Гельдера з показником $p$ та означенням норми \eqref{Brpt-period-dek1-bs}, продовжимо оцінку \eqref{T3-J4-os1}
$$
J_4 = \sum\limits_{(\boldsymbol{s},\boldsymbol{\gamma}) \geqslant n} 2^{-\left((\boldsymbol{s},\boldsymbol{r})- \|s\|_1\left(\frac{1}{p}-\frac{1}{q}\right)\right)} 2^{(\boldsymbol{s},\boldsymbol{r})} \|\delta_{\boldsymbol{s}}(f)\|_p \leqslant
$$
$$
\leqslant \left(\sum\limits_{(\boldsymbol{s},\boldsymbol{\gamma}) \geqslant n} 2^{(\boldsymbol{s},\boldsymbol{r})p}\|\delta_{\boldsymbol{s}}(f)\|^p_p\right)^{\frac{1}{p}} \left(\sum\limits_{(\boldsymbol{s},\boldsymbol{\gamma}) \geqslant n} 2^{-\left((\boldsymbol{s},\boldsymbol{r})- \|s\|_1\left(\frac{1}{p}-\frac{1}{q}\right)\right)p'}\right)^{\frac{1}{p'}} \ll
 $$
\begin{equation}\label{T3-J4-os2}
 \ll \|f\|_{B^{\boldsymbol{r}}_{p,p}} \left(\sum\limits_{(\boldsymbol{s},\boldsymbol{\gamma})\geqslant n} 2^{-\left(\boldsymbol{s},\boldsymbol{r}-\frac{1}{p}+\frac{1}{q}\right)p'}\right)^{\frac{1}{p'}} \leqslant \left(\sum\limits_{(\boldsymbol{s},\boldsymbol{\gamma}) \geqslant n} 2^{-(\boldsymbol{s},\boldsymbol{\widetilde{\gamma}})\left(r_1-\frac{1}{p}+\frac{1}{q}\right)p'}\right)^{\frac{1}{p'}},
\end{equation}
де $\frac{1}{p}+\frac{1}{p'}= 1$, $\boldsymbol{r}-\frac{1}{p}+\frac{1}{q}$~--- вектор з координатами $r_j-\frac{1}{p}+\frac{1}{q}$, $j= \overline{1,d}$, і  відповідно $\boldsymbol{\widetilde{\gamma}}=(\widetilde{\gamma}_1,\ldots,\widetilde{\gamma}_d)$, $\widetilde{\gamma}_j={\displaystyle {\frac{r_j-\frac{1}{p}+\frac{1}{q}}{r_1-\frac{1}{p}+\frac{1}{q}}}}$, $j=\overline{1,d}$. Легко бачити, що  $\widetilde{\gamma}_j = \gamma_j$, $j=\overline{1,\nu}$ і $1<\gamma_j<\widetilde{\gamma}_j$, $j=\overline{\nu + 1, d}$. Враховуючи це і скориставшись лемою~А, по відношенню до останньої суми в \eqref{T3-J4-os2},  отримаємо
\begin{equation}\label{T3-J4-os3}
J_4 \ll  2^{-n\left(r_1-\frac{1}{p}+\frac{1}{q}\right)} n^{(\nu - 1)\left(1 - \frac{1}{p}\right)}.
\end{equation}

Насамкінець співставляючи \eqref{T3-J4-os1}--\eqref{T3-J4-os3}, отримаємо  співвідношення
$$
E_{Q^{\boldsymbol{\gamma}}_n}\big(W^{\boldsymbol{r}}_{p,\boldsymbol{\alpha}}\big)_{B_{q,1}} \leqslant \mathscr{E}_{Q^{\boldsymbol{\gamma}}_n}\big(W^{\boldsymbol{r}}_{p,\boldsymbol{\alpha}}\big)_{B_{q,1}} \ll \mathscr{E}_{Q^{\boldsymbol{\gamma}}_n}\big(B^{\boldsymbol{r}}_{p,p}\big)_{B_{q,1}} \ll 2^{-n\left(r_1-\frac{1}{p}+\frac{1}{q}\right)} n^{(\nu - 1)\left(1-\frac{1}{p}\right)}.
$$

Переходячи в \eqref{T3-En-Wp-Bq1} до встановлення оцінок знизу зауважимо, що їх, як і у теоремі~1, достатньо отримати у випадку $\boldsymbol{\alpha}=\boldsymbol{0}$ і $\nu = d$. Для цього також будемо розглядати  функцію $g_1$, яка, як показано при доведенні теореми~1, належить класу $W^{\boldsymbol{r}}_p$, ${\boldsymbol{r}=(r_1,\ldots,r_1) \in \mathbb{R}^d}$.

Врахувавши, що $S_{Q^{\boldsymbol{1}}_n}(g_1) = 0$ будемо мати
$$
E_{Q^{\boldsymbol{1}}_n}\big(g_1\big)_{B_{q,1}} \asymp \mathscr{E}_{Q^{\boldsymbol{1}}_n}\big(g_1\big)_{B_{q,1}} = \|g_1\|_{B_{q,1}}\asymp 2^{-n \left(r_1+1-\frac{1}{p}\right)} n^{-\frac{d-1}{p}} \sum\limits_{(\boldsymbol{s},\boldsymbol{1})=n} \|\delta_{\boldsymbol{s}}(d_n)\|_q \asymp
$$
$$
\asymp 2^{-n\left(r_1+1-\frac{1}{p}\right)} n^{-\frac{d - 1}{p}} \sum\limits_{(\boldsymbol{s},\boldsymbol{1})=n} 2^{(\boldsymbol{s},\boldsymbol{1})\left(1-\frac{1}{q}\right)} \asymp
$$
$$
\asymp 2^{-n\left(r_1+1-\frac{1}{p}\right)} n^{-\frac{d - 1}{p}} 2^{n\left(1-\frac{1}{q}\right)} n^{d-1} =  2^{-n\left(r_1-\frac{1}{p}+\frac{1}{q}\right)} n^{(d - 1)\left(1 - \frac{1}{p}\right)}.
$$

Оцінки знизу встановлено.

Теорему~3 доведено. \vskip 2 mm

Наведемо відповідне до теореми~3 твердження у просторі $L_q$.

\textbf{Теорема~Д} \cite{Galeev-1978}\textbf{.} \it Нехай $d \geqslant 2$, $1 < p < q < \infty$, $r_1 > \frac{1}{p}-\frac{1}{q}$. Тоді при $\boldsymbol{\alpha} \in \mathbb{R}^d$ справедливі співвідношення
\begin{equation}\label{T3-En-Wp-Lq}
E_{Q^{\boldsymbol{\gamma}}_n}\big(W^{\boldsymbol{r}}_{p,\boldsymbol{\alpha}}\big)_q \asymp \mathscr{E}_{Q^{\boldsymbol{\gamma}}_n}\big(W^{\boldsymbol{r}}_{p,\boldsymbol{\alpha}}\big)_q
\asymp 2^{-n\left(r_1-\frac{1}{p}+\frac{1}{q}\right)}.
\end{equation} \rm

Отже, співставивши оцінки \eqref{T3-En-Wp-Bq1} і \eqref{T3-En-Wp-Lq} при $2 \leqslant p < q <\infty$ виявляємо, що відповідні апроксимаційні характеристик класів $W^{\boldsymbol{r}}_{p,\boldsymbol{\alpha}}$ у просторі
$L_p$ за винятком випадку $\nu = 1$ відрізняються за порядком. Крім цього у просторі $B_{q, 1}$ прослідковується залежність одержаних оцінок від показника розмірності $\nu$. Інша ситуація спостерігається в одновимірному випадку.

Справедливим є таке твердження.

\textbf{Теорема~3$'$.} \it Нехай $d =1$, $1<p<q<\infty$, $r >\frac{1}{p}-\frac{1}{q}$, $\alpha \in \mathbb{R}$. Тоді виконуються співвідношення
\begin{equation}\label{T3-En-Wp-Bq1-d1}
E_{2^n}\big(W^r_{p,\alpha}\big)_{B_{q,1}} \asymp \mathscr{E}_{2^n}\big(W^r_{p,\alpha})_{B_{q,1}} \asymp 2^{-n\left(r-\frac{1}{p}+\frac{1}{q}\right)}.
\end{equation} \rm

\emph{\textbf{Доведення.}} Оцінки  зверху  в \eqref{T3-En-Wp-Bq1-d1} випливають із оцінок відповідних характеристик    класів $H^r_p$ (див. Наслідок~1 з \cite{RomanyukAV-Pozharska-2023-CMP}) і вкладення $W^{r}_{p,\alpha} \subset H^r_p$.

Оцінка знизу для найкращого наближення $E_{2^n}\big(W^{r}_{p,\alpha}\big)_{B_{q,1}}$ є наслідком оцінки цієї величини у просторі $L_q$ (див. \cite[гл.\,1, \S\,3]{Temlyakov-1993m}) і співвідношення $\|\cdot\|_{B_{q,1}} \gg \|\cdot\|_q$.

Теорему~3$'$ доведено. \vskip 2 mm

Зауважимо, з наведених міркувань легко бачити, що розглянуті характеристики класів $W^r_{p,\alpha}$ у просторах $B_{q,1}$ і $L_q$ однакові за порядком.

На завершення розглянемо ще одне співвідношення між параметрами $p$ і $q$.

\textbf{Теорема~4.} \it Нехай $d \geqslant 2$, $1\leqslant q \leqslant 2$, $q < p < \infty$, $r_1 >0$. Тоді при $\boldsymbol{\alpha}\in \mathbb{R}^d$ справедливі співвідношення
\begin{equation}\label{T4-En-Wp-Bq1}
E_{Q^{\boldsymbol{\gamma'}}_n}\big(W^{\boldsymbol{r}}_{p,\boldsymbol{\alpha}}\big)_{B_{q,1}} \asymp   \mathscr{E}_{Q^{\boldsymbol{\gamma'}}_n}\big(W^{\boldsymbol{r}}_{p,\boldsymbol{\alpha}}\big)_{B_{q,1}}
 \asymp 2^{-nr_1} n^{\frac{\nu-1}{2}}.
\end{equation} \rm

\emph{\textbf{Доведення.}} При встановленні оцінок зверху в \eqref{T4-En-Wp-Bq1} розглянемо два випадки.

а) Нехай $1\leqslant q <p\leqslant2$. Тоді оцінки зверху величини $\mathscr{E}_{Q^{\boldsymbol{\gamma'}}_n}\big(W^{\boldsymbol{r}}_{p,\boldsymbol{\alpha}}\big)_{B_{q,1}}$ є наслідком теореми~1, оскільки
$$
\mathscr{E}_{Q^{\boldsymbol{\gamma'}}_n}\big(W^{\boldsymbol{r}}_{p,\boldsymbol{\alpha}}\big)_{B_{q,1}} \leqslant \mathscr{E}_{Q^{\boldsymbol{\gamma'}}_n}\big(W^{\boldsymbol{r}}_{p,\boldsymbol{\alpha}}\big)_{B_{p,1}}
 \asymp 2^{-nr_1} n^{\frac{\nu-1}{2}}.
$$

б) Нехай $1\leqslant q \leqslant2<p<\infty$. Тоді врахувавши, що $W^{\boldsymbol{r}}_{p,\boldsymbol{\alpha}}\subset W^{\boldsymbol{r}}_{2,\boldsymbol{\alpha}}$ і знову скориставшись результатом теореми~1 будемо мати
$$
\mathscr{E}_{Q^{\boldsymbol{\gamma'}}_n}\big(W^{\boldsymbol{r}}_{p,\boldsymbol{\alpha}}\big)_{B_{q,1}} \leqslant \mathscr{E}_{Q^{\boldsymbol{\gamma'}}_n}\big(W^{\boldsymbol{r}}_{2,\boldsymbol{\alpha}}\big)_{B_{2,1}}
 \asymp 2^{-nr_1} n^{\frac{\nu-1}{2}}.
$$

Оцінки зверху встановлено.

Оцінки знизу в \eqref{T4-En-Wp-Bq1} отримуються як наслідок оцінки колмогоровського поперечника $d_M\big(W^{\boldsymbol{r}}_{p,\boldsymbol{\alpha}}, B_{1,1}\big)$ (Теорема~А), тобто при $M\asymp 2^n n^{\nu-1}$ матимемо
$$
\mathscr{E}_{Q^{\boldsymbol{\gamma'}}_n}\big(W^{\boldsymbol{r}}_{p,\boldsymbol{\alpha}}\big)_{B_{q,1}} \geqslant E_{Q^{\boldsymbol{\gamma'}}_n}\big(W^{\boldsymbol{r}}_{p,\boldsymbol{\alpha}}\big)_{B_{q,1}} \gg d_M\big(W^{\boldsymbol{r}}_{p,\boldsymbol{\alpha}}, B_{q,1}\big) \geqslant d_M\big(W^{\boldsymbol{r}}_{p,\boldsymbol{\alpha}}, B_{1,1}\big) \asymp
$$
$$
\asymp M^{- r_1} (\log^{\nu -1 }M)^{r_1 + \frac{1}{2}} \asymp 2^{-nr_1}n^{\frac{\nu-1}{2}}.
$$

Теорему~4 доведено. \vskip 2 mm

Твердження, відповідне теоремі~4 у просторі $L_q$, має такий вигляд.

\textbf{Теорема~Е} \cite[гл.\,3, \S\,3]{Temlyakov-1993m}\textbf{.} \it Нехай $d \geqslant 2$, $1<q < p < \infty$, $r_1 >0$. Тоді при $\boldsymbol{\alpha}\in \mathbb{R}^d$ виконуються співвідношення
\begin{equation}\label{T4-En-Wp-Lq}
E_{Q^{\boldsymbol{\gamma'}}_n}\big(W^{\boldsymbol{r}}_{p,\boldsymbol{\alpha}}\big)_{q} \asymp   \mathscr{E}_{Q^{\boldsymbol{\gamma'}}_n}\big(W^{\boldsymbol{r}}_{p,\boldsymbol{\alpha}}\big)_{q}
 \asymp 2^{-nr_1}.
\end{equation} \rm

Таким чином, порівнявши співвідношення \eqref{T4-En-Wp-Bq1} і \eqref{T4-En-Wp-Lq}, бачимо, що порядкові оцінки відповідних величин при $\nu \neq 1$ є різними.

У доповнення до теореми~4 наведемо результат стосовно одновимірного випадку.

\textbf{Теорема~4$'$.} \it Нехай $d =1$, $1 <q<p<\infty$, $r_1 >0$. Тоді при $\alpha \in \mathbb{R}$
$$
E_{2^n}\big(W^r_{p,\alpha}\big)_{B_{q,1}} \asymp \mathscr{E}_{2^n}\big(W^r_{p,\alpha}\big)_{B_{q,1}}
  \asymp  2^{ -n  r_1}.
$$\rm

\emph{\textbf{Доведення.}} Оцінка зверху величини $\mathscr{E}_{2^n}\big(W^r_{p,\alpha}\big)_{B_{q,1}}$ випливає із співвідношення
$$
\mathscr{E}_{2^n}(H^{r_1}_{p})_{B_{q,1}} \asymp 2^{ -n  r_1}
$$
(див. Наслідок~1 \cite{RomanyukAV-Pozharska-2023-CMP}) згідно із вкладенням $W^r_{p,\alpha} \subset H^r_p$.

Оцінка знизу величини  $E_{2^n}\big(W^r_{p,\alpha}\big)_{B_{q,1}}$ є наслідком співвідношення
$$
E_{2^n}\big(W^r_{p,\alpha}\big)_q \asymp 2^{-nr}
$$
(див. \cite[Гл.\,1, \S\,3]{Temlyakov-1993m}) і нерівності $\|\cdot\|_{B_{q,1}} \gg \|\cdot\|_q$.

Теорему~4$'$ доведено. \vskip 2 mm

Отже, в одновимірному випадку при $1 < q < p < \infty$ відповідні апроксимаційні характеристики класів $W^r_{p,\alpha}$ у просторах $B_{q,1}$ і $L_q$ є однаковими за порядком.

\vskip 5 mm

\textbf{Acknowledgments:} \emph{This work was supported by a grant from the Simons Foundation (1290607,  A.\,S.~Romanyuk, S.\,Ya.~Yanchenko).}

\vskip 3 mm

\textbf{Contact information:}
Department of the Theory of Functions of the Institute of Mathematics of the National
Academy of Sciences of Ukraine, 3, Tereshenkivska st., 01024, Kyiv, Ukraine.

\vskip 3 mm


\vskip 2 mm

\begin{thebibliography}{10}

\bibitem{Temlyakov-1990Proc-189}   V.\,N.~Temlyakov, {\it Estimates of the asymptotic  characteristics of classes  of  functions  with  bounded mixed derivative or difference}, Proc. Steklov Inst. Math., {\bf 189},  161\,--\,197 (1990).

\bibitem{Kashin-Temlyakov-1994MN} B.\,S.~Kashin, V.\,N.~Temlyakov,  {\it On best $m$-term approximations and the entropy of sets in the space $L_1$}, Math. Notes, \textbf{56}, No.\,5, 1137\,--\,1157 (1994). https://doi.org/10.1007/BF02274662

\bibitem{Belinskii-1998JAT}  E.\,S.~Belinsky, {\it Estimates of entropy numbers and Gaussian measures for classes of functions with bounded mixed derivative}, J. Approx. Theory, \textbf{93}, No.\,1, 114\,--\,127 (1998) https://doi.org/10.1006/jath.1997.3157

\bibitem{RomanyukAV-19UMJ2} A.\,S.~Romanyuk, V.\,S.~Romanyuk, {\it Approximating characteristics of the classes of periodic multivariate functions in the space $B_{\infty,1}$}, Ukrainian Math. J., \textbf{71}, No.\,2, 308\,--\,321 (2019). https://doi.org/10.1007/s11253-019-01646-3

\bibitem{RomanyukAV-19UMJ8} A.\,S.~Romanyuk, V.\,S.~Romanyuk, {\it Estimation of some approximating characteristics of the classes of periodic functions of one and many variables}, Ukrainian Math. J., \textbf{71}, No.\,8, 1257\,--\,1272 (2020).  https://doi.org/10.1007/s11253-019-01711-x

\bibitem{RomanyukAV-21UMB} A.\,S.~Romanyuk, V.\,S.~Romanyuk, {\it Approximative characteristics and properties of operators of the best approximation of classes of functions from the Sobolev and Nikol’skii--Besov spaces}, J. Math. Sci., \textbf{252}, No.\,4, 508\,--\,525 (2021). https://doi.org/10.1007/s10958-020-05177-2

\bibitem{Hembarska-Zaderei-2022} S.\,B.~Hembars’ka,  P.\,V.~Zaderei, {\it Best orthogonal trigonometric approximations of the Nikol’skii--Besov type classes of periodic functions in the fpace $B_{\infty,1}$},  Ukrainian Math. J., \textbf{74}, No.\,6, 883\,--\,895 (2022). https://doi.org/10.1007/s11253-022-02115-0

\bibitem{Romanyuk-Yanchenko-2021UMJ} A.\,S.~Romanyuk, S.\,Ya.~Yanchenko, {\it  Estimates of approximating characteristics and the properties of the operators of best approximation for the classes of periodic functions in the space $B_{1,1}$}, Ukrainian Math. J., \textbf{73}, No.\,8, 1278\,--\,1298 (2022). https://doi.org/10.1007/s11253-022-01990-x

\bibitem{Romanyuk-Yanchenko-2022-6UMJ} A.\,S.~Romanyuk, S.\,Ya.~Yanchenko, {\it Approximation of the classes of periodic functions of one and many variables from the Nikol’skii--Besov and Sobolev spaces},  Ukrainian Math. J., \textbf{74}, No.\,6, 967\,--\,980 (2022). https://doi.org/10.1007/s11253-022-02110-5

\bibitem{Lizorkin-Nikolsky-1989} P.\,I.~Lizorkin, S.\,M.~Nikol’skii,  {\it Function spaces of mixed smoothness from the decomposition point of view}, Proc. Steklov Inst. Math., \textbf{187}, 163\,--\,184 (1990).

\bibitem{Temlyakov-1986m}  V.\,N.~Temlyakov,  {\it Approximation of functions with bounded mixed derivative},  Proc. Steklov Inst. Math.,  \textbf{178}, 1\,--\,121 (1989).

\bibitem{Temlyakov-1993m} V.\,N.~Temlyakov, {\it Approximation of periodic functions}, New York: Nova Sci. Publ. Inc. (1993).

\bibitem{Trigub-Belinsky} R.\,M.~Trigub, E.\,S.~Belinsky,  {\it Fourier Analysis and Approximation of Functions}, Kluwer Academic Publishers, Dordrecht, 2004. https://doi.org/10.1007/978-1-4020-2876-2

\bibitem{Romanyuk-2012m} A.\,S.~Romanyuk, {\it Approximating Characteristics of the Classes of Periodic Functions of Many Variables}, Proc. of the Institute of Mathematics, National Academy of Sciences of Ukraine, Kyiv, 2012.

\bibitem{Temlyakov-2018} V.\,N.~Temlyakov, {\it Multivariate approximation}, Cambridge University Press, Cambridge, 2018. https://doi.org/10.1017/9781108689687


\bibitem{Cross-2018} D.~D\~{u}ng, V.~Temlyakov and T.~Ullrich, {\it Hyperbolic Cross Approximation},  Advanced Courses in Mathematics, CRM Barselona, Birkh\"{a}user/Springer, Cham, 2018.  https://doi.org/10.1007/978-3-319-92240-9

\bibitem{Fedunyk-Hembarska-2020-CMP} O.\,V.~Fedunyk-Yaremchuk, M.\,V.~Hembars’kyi, S.\,B.~Hembars’ka, {\it Approximative characteristics of the Nikol’skii--Besov-type classes of periodic functions in the space $B_{\infty,1}$}, Carpathian Math. Publ.,  \textbf{12}, No.\,2, 376\,--\,391 (2020). https://doi.org/10.15330/cmp.12.2.376-391

\bibitem{RomanyukAV-Pozharska-2023-CMP}  A.\,S.~Romanyuk, V.\,S.~Romanyuk, K.\,V.~Pozharska,  S.\,B.~Hembars'ka,  {\it Characteristics of linear and nonlinear approximation of isotropic classes of periodic multivariate functions}, Carpathian Math. Publ.,  \textbf{15}, No.\,1, 78\,--\,94 (2023). https://doi.org/10.15330/cmp.15.1.78-94

\bibitem{Hembarska-Fedunyk-2023} S.\,B.~Hembars’ka, I.\,A.~Romanyuk,  O.\,V.~Fedunyk-Yaremchuk,  {\it Characteristics of the linear and nonlinear approximations of the Nikol’skii--Besov-type classes of periodic functions of several variables}, J. Math. Sci., \textbf{274}, No.\,3, 307\,--\,326 (2023). https://doi.org/10.1007/s10958-023-06602-y

\bibitem{Kolmogorof_1936} A.~Kolmogoroff, {\it \"{U}ber die beste Ann\"{a}herung von Functionen einer gegeben Functionenclasse}, Ann. Math., \textbf{37}, No.\,1, 107\,--\,110 (1936). https://doi.org/10.2307/1968691

\bibitem{Romanyuk-UMJ17-10} A.\,S.~Romanyuk,  {\it Entropy numbers and widths for the classes $B^r_{p,\theta}$
 of periodic functions of many variables}, Ukrainian Math. J., \textbf{68}, No.\,10, 1620\,--\,1636 (2017). https://doi.org/10.1007/s11253-017-1315-9

\bibitem{Nikolsky_51} S.\,M.~Nikol'skii, { \it Inequalities for entire functions of finite power and their application to the theory of differentiable functions of many variables}, Tr. Mat. Inst. Steklova, \textbf{38}, 244\,--\,278 (1951).

\bibitem{Mityagin-1962}    B.\,S.~Mityagin,  {\it The approximation   of  functions   in the    space $L^p$  and  $C$  on  the   torus},   Mat. Sb., \textbf{58}, No.\,4, 397\,--\,414 (1962).

\bibitem{Nikolskaya-75-sib} N.\,S.~Nikol'skaya,  {\it The approximation in the $L_p$ metric of differentiable functions of several variables by Fourier sums}, Sib. Math. J., \textbf{15}, No.\,2, 282\,--\,295 (1974). https://doi.org/10.1007/BF00968291

\bibitem{Galeev-1978} E.\,M.~Galeev, {\it Approximation by Fourier sums of classes of functions with several bounded derivatives}, Math. Notes, \textbf{23}, No.\,2, 109\,--\,117 (1978). https://doi.org/10.1007/BF01153149

\end{thebibliography}
\end{document}